%% file: cube.tex
\documentclass{amsart}
\usepackage{amssymb}

\usepackage{psfig}

\theoremstyle{plain}
\newtheorem{theorem}{Theorem}[section]
\newtheorem{corollary}[theorem]{Corollary}
\newtheorem{lemma}[theorem]{Lemma}
\newtheorem{proposition}[theorem]{Proposition}
\theoremstyle{definition}
\newtheorem{conjecture}[theorem]{Conjecture}
\theoremstyle{remark}
\newtheorem{acknowledgements}{Acknowledgements}

\numberwithin{equation}{section}

\newcommand{\figlabone}{\relax}
\newcommand{\figlabtwo}{\relax}
\input inner

\input tcibites

\begin{document}
\title[Spectral pairs in Cartesian coordinates]{Spectral pairs in Cartesian coordinates}
\author{Palle E. T. Jorgensen}
\address{Department of Mathematics\\
The University of Iowa\\
Iowa City, IA 52242\\
U.S.A.}
\email{jorgen@math.uiowa.edu}
\author{Steen Pedersen}
\address{Department of Mathematics\\
Wright State University\\
Dayton, OH 45435\\
U.S.A.}
\email{steen@math.wright.edu}
\thanks{Work supported by the National Science Foundation.}
\dedicatory{Dedicated to the memory of Irving E. Segal}
\subjclass{42C05, 22D25, 46L55, 47C05}
\keywords{Spectral pair, translations, tilings, Fourier basis, operator extensions,
induced representations, spectral resolution, Hilbert space.}

\begin{abstract}
Let $\Omega \subset \mathbb{R}^{d}$ have finite positive Lebesgue measure,
and let $\mathcal{L}^{2}\left( \Omega \right) $ be the corresponding Hilbert
space of $\mathcal{L}^{2}$-functions on $\Omega $. We shall consider the
exponential functions $e_{\lambda }$ on $\Omega $ given by
$e_{\lambda }\left( x\right)  =e^{i2\pi \lambda \cdot x}$.
If these functions form an orthogonal basis for $\mathcal{L}^{2}\left(
\Omega \right) $, when $\lambda $ ranges over some subset $\Lambda $ in $%
\mathbb{R}^{d}$, then we say that $\left( \Omega ,\Lambda \right) $ is a 
\emph{spectral pair}, and that $\Lambda $ is a \emph{spectrum.} 
We conjecture that $\left( \Omega ,\Lambda \right) $ is a spectral pair
if and only if the translates of some set $\Omega ^{\prime}$
by the vectors of $\Lambda $ tile $\mathbb{R}^{d}$.
In the special case of 
$\Omega =I^{d}$, the $d$-dimensional unit cube,
we prove this conjecture,
with $\Omega ^{\prime}=I^{d}$,
for $d\leq 3$,
describing all the tilings by $I^{d}$,
and for all
$d$ when $\Lambda $ is a discrete
periodic set.
In an appendix we generalize the notion
of spectral pair to measures on a locally
compact abelian group and its dual.
\end{abstract}

\maketitle\renewcommand{\figlabone}{\setlength{\unitlength}{54pt}\smash{\makebox[0pt]{\begin{picture}(4,5)(-2.71,0.02)\put(-1,-0.1){\makebox(0.25,1)[lt]{$\beta _{1}-\beta _{0}$}}\put(0,0.11){\makebox(0.25,1)[lt]{$\beta _{2}-\beta _{1}$}}\put(1,0.36){\makebox(0.25,1)[lt]{$\beta _{3}-\beta _{2}$}}\end{picture}}}}\renewcommand{\figlabtwo}{\setlength{\unitlength}{54pt}\smash{\makebox[0pt]{\begin{picture}(5,4)(-0.65,-2.04)\put(-0.1,-1){\makebox(0.25,1)[rt]{$\beta _{1}-\beta _{0}$}}\put(0.11,0){\makebox(0.25,1)[rt]{$\beta _{2}-\beta _{1}$}}\put(0.36,1){\makebox(0.25,1)[rt]{$\beta _{3}-\beta _{2}$}}\end{picture}}}}%

\section{\label{S1}Introduction}

The setting of \emph{spectral pairs} in $d$ real dimensions involves two
subsets $\Omega $ and $\Lambda $ in $\mathbb{R}^{d}$ such that $\Omega $ has
finite and positive $d$-dimensional Lebesgue measure, and $\Lambda $ is an
index set for an orthogonal $\mathcal{L}^{2}\left( \Omega \right) $-basis $%
e_{\lambda }$ of exponentials, i.e., 
\begin{equation}
e_{\lambda }\left( x\right) =e^{i2\pi \lambda \cdot x},\quad x\in \Omega
,\;\lambda \in \Lambda  \label{eq1}
\end{equation}
where $\lambda \cdot x=\sum_{j=1}^{d}\lambda _{j}x_{j}$. We use vector
notation $x=\left( x_{1},\cdots ,x_{d}\right) $, $\lambda =\left( \lambda
_{1},\cdots ,\lambda _{d}\right) $, $x_{j},\lambda _{j}\in \mathbb{R}$, $%
j=1,\dots ,d$. The basis property refers to the Hilbert space $\mathcal{L}%
^{2}\left( \Omega \right) $ with inner product 
\begin{equation}
\ip{f}{g}%
_{\Omega }:=\int_{\Omega }\overline{f\left( x\right) }g\left( x\right) \,dx
\label{eq2}
\end{equation}
where $dx=dx_{1}\cdots dx_{d}$, and $f,g\in \mathcal{L}^{2}\left( \Omega
\right) $. The corresponding norm is 
\begin{equation}
\left\| f\right\| _{\Omega }^{2}:=%
\ip{f}{f}%
_{\Omega }=\int_{\Omega }\left| f\left( x\right) \right| ^{2}\,dx,
\label{eq3}
\end{equation}
as usual. It follows that the spectral pair property for a pair $\left(
\Omega ,\Lambda \right) $ is equivalent to 
the nonzero elements of the set
\begin{equation*}
\Lambda -\Lambda =\left\{ \lambda -\lambda ^{\prime }:\lambda ,\lambda
^{\prime }\in \Lambda \right\}
\end{equation*}
being contained in the \emph{zero-set} of the complex valued function 
\begin{equation}
z\longmapsto \int_{\Omega }e^{i2\pi z\cdot x}\,dx=:F_{\Omega }\left( z\right)
\label{eq4}
\end{equation}
and the corresponding $e_{\lambda }$-set $%
\left\{ e_{\lambda }:\lambda \in \Lambda \right\} $ being \emph{total} in $%
\mathcal{L}^{2}\left( \Omega \right) $. Recall, totality means that the span
of the $e_{\lambda }$'s is dense in $\mathcal{L}^{2}\left( \Omega \right) $
relative to the $\left\| \,\cdot \,\right\| _{\Omega }$-norm, or,
equivalently, that $f=0$ is the only $\mathcal{L}^{2}\left( \Omega \right) $%
-solution to: 
\begin{equation*}
\ip{f}{e_\lambda }%
_{\Omega }=0\text{,\quad for all }\lambda \in \Lambda .
\end{equation*}
Note, $F_{\Omega }\left( z\right) $ is defined for any $z=\left(
z_{1},\cdots ,z_{d}\right) \in \mathbb{R}^{d}$ since $\Omega $ has finite
measure and $e^{i2\pi z\cdot x}$ has absolute value $=1$.
We refer to the
book \cite{Ped97} for a summary of the theory of \emph{spectral pairs.} It
was developed in the previous joint papers by the coauthors
\cite{JoPe87,JoPe91,JoPe92,JoPe93a,JoPe93b,JoPe94,JoPe95,JoPe96}
and elsewhere, e.g.,
\cite{LRW98,LaWa96c,LaWa97a,LaWa97b}.
We recall that
Fuglede showed \cite{Fug74} that the disk and the triangle in two dimensions
are \emph{not spectral sets,} in the sense that, if $\Omega $ is one of
these sets, then there is \emph{no} possible choice for $\Lambda $ such that 
$\left( \Omega ,\Lambda \right) $ is a spectral pair in $\mathbb{R}^{2}$.

If $\Omega \subset \mathbb{R}^{d}$ is open, then we consider the partial
derivatives $\frac{\partial \;}{\partial x_{j}}$, $j=1,\dots ,d$, defined on 
$C_{c}^{\infty }\left( \Omega \right) $ as unbounded skew-symmetric
operators in $\mathcal{L}^{2}\left( \Omega \right) $. The corresponding
versions $\frac{1}{\sqrt{-1}}\frac{\partial \;}{\partial x_{j}}$ are
symmetric of course. We say that $\Omega $ has the \emph{extension property}
if there are commuting self-adjoint extension operators $H_{j}$, i.e., 
\begin{equation}
\frac{1}{i}\frac{\partial \;}{\partial x_{j}}\subset H_{j},\quad j=1,\dots
,d.  \label{eq7}
\end{equation}

We have (see \cite{Fug74,Jor82,Ped87,JoPe92})

\begin{theorem}
\label{Thm1.1}\textup{(Fuglede, Jorgensen, Pedersen)} Let $\Omega \subset %
\mathbb{R}^{d}$ be open and connected with finite and positive Lebesgue
measure. Then $\Omega $ has the extension property if and only if it is a
spectral set. If $\Omega $ is only assumed open, then the spectral-set
property implies the extension property, but not conversely.
\end{theorem}


Some of the interest in spectral pairs derives from their connection to 
\emph{tilings.} A subset $\Omega \subset \mathbb{R}^{d}$ with nonzero
measure is said to be a \emph{tile%
}
if there is a set $%
L\subset \mathbb{R}^{d}$ such that the translates $\left\{ \Omega +l:l\in
L\right\} $ cover $\mathbb{R}^{d}$ up to measure zero, and if the
intersections 
\begin{equation}
\left( \Omega +l\right) \cap \left( \Omega +l^{\prime }\right) \text{\quad
for }l\ne l^{\prime }\text{ in }L  \label{eq10}
\end{equation}
have measure zero. We will call $\left( \Omega ,L\right) $ a \emph{tiling
pair} and we will say that $L$ is a \emph{tiling set}. The \emph{%
Spectral-Set conjecture} due to Fuglede (see \cite
{Fug74,Jor82,Ped87,JoPe95,LaWa96c,LaWa97a,LaWa97b}) states:

\begin{conjecture}
\label{Con1.3}Let $\Omega \subset \mathbb{R}^{d}$ have positive and finite
Lebesgue measure. Then $\Omega $ is a spectral set if and only if $\Omega $
is a tile, i.e., there exists a set $L$ so that $(\Omega ,L)$ is a spectral
pair if and only if there exists a set $L^{\prime }$ so that $(\Omega %
,L^{\prime })$ is a tiling pair.
\end{conjecture}

We formulate a ``dual'' conjecture.

\begin{conjecture}
\label{Con1.4}Let $L\subset \mathbb{R}^{d}$. Then $L$ is a spectrum if and
only if $L$ is a tiling set, i.e., there exists a set $\Omega $ so that $(%
\Omega ,L)$ is a spectral pair if and only if there exists a set $\Omega %
^{\prime }$ so that $(\Omega ^{\prime },L)$ is a tiling pair.
\end{conjecture}

\begin{conjecture}
\label{Con1.5}Let $L\subset \mathbb{R}^{d}$. Then $\left( I^{d},L\right) $
is a spectral pair if and only if $\left( I^{d},L\right) $ is a tiling pair.
\end{conjecture}

The significance of the special case $\Omega =I^{d}$ lies in part in the
results below where we show, for $d=1,2,3$, that $\left( I^{d},\Lambda
\right) $, $\Lambda \subset \mathbb{R}^{d}$, is a \emph{spectral pair} if
and only if $I^{d}$ \emph{tiles }$\mathbb{R}^{d}$\emph{\ by }$\Lambda
$\emph{-translates}. Our proofs
also
construct all possible spectra for the unit cube
when $d=1,2,3$. In Section \ref{S:periodic} we
establish Conjecture \ref{Con1.5} for all $d$ when
$\Lambda $ is a discrete periodic set.

Tiling questions for $I\subset \mathbb{R}$ are trivial, but not so for $%
I^{d}\subset \mathbb{R}^{d}$ when $d\ge 2$. The connection between \emph{%
tiles} and \emph{spectrum} is more direct for $\Omega =I^{d}$ than for other
examples of sets $\Omega $.
This is
explained by the following (easy) lemma relating the problems to the
function $F_{\Omega }$ from (\ref{eq4}) above.

\begin{lemma}
\label{Lem1.4}If $\Omega =I^{d}$, then the zero-set for the function $F_{%
\Omega }$ in \textup{(\ref{eq4})} is 
\begin{equation}
\mathbf{Z}_{I^{d}}=\left\{ z\in \mathbb{C}^{d}\diagdown \left\{ 0\right\}
:\exists j\in \left\{ 1,\dots ,d\right\} \;\mathrm{s.t.}\;z_{j}\in \mathbb{Z}%
\diagdown \left\{ 0\right\} \right\} .  \label{eq11}
\end{equation}
\end{lemma}

\begin{proof}%
The function $F_{I^{d}}\left( \,\cdot \,\right) $ factors as follows. 
\begin{equation}
F_{I^{d}}\left( z\right) =\prod_{j=1}^{d}\frac{e^{i2\pi z_{j}}-1}{i2\pi z_{j}%
}  \label{eq12}
\end{equation}
for $z=\left( z_{1},\dots ,z_{d}\right) \in \mathbb{C}^{d}$, with the
interpretation that the function $z\mapsto \frac{e^{i2\pi z}-1}{i2\pi z}$ is 
$1$ when $z=0$ in $\mathbb{C}$.%
\end{proof}%

In particular, if $(I^{d},\Lambda )$ is a spectral pair, then $\Lambda
-\Lambda \subset \mathbf{Z}_{I^{d}}\cup \{0\}$.
The corresponding result for
tilings is non-trivial, it was proved by Keller \cite{Kel30,Kel37}, a
detailed proof appears in \cite{Per40}. The precise statement of Keller's
theorem is:

\begin{theorem}
\label{T:Keller}If $(I^{d},\Lambda )$ is a tiling pair, then $\Lambda
-\Lambda \subset \mathbf{Z}_{I^{d}}\cup \{0\}$, where $\mathbf{Z}_{I^{d}}$
is given by \textup{(\ref{eq11})}.
\end{theorem}

Let $\mu ,\nu $ be two Borel measures on $\mathbb{R}^{d}$. We will say that $%
(\mu ,\nu )$ is a \emph{tiling pair} if the convolution, $\mu *\nu $, of $%
\mu $ and $\nu $ is Lebesgue measure on $\mathbb{R}^{d}$. This coincides
with the previous definition of a tiling pair in the sense that if $(\Omega %
,L)$ is a pair of subsets of $\mathbb{R}^{d}$ so that $\Omega $ has finite
positive Lebesgue measure, $L$ is discrete, $\omega $ denotes Lebesgue
measure restricted to $\Omega $, and $\ell $ denotes counting measure on $L$%
, then $(\Omega ,L)$ is tiling pair if and only if $(\omega ,\ell )$ is a
tiling pair. Since convolution is commutative, $(\mu ,\nu )$ is a tiling
pair if and only if $(\nu ,\mu )$ is a tiling pair. In the appendix we
introduce (and investigate properties of) a notion of a spectral pair of
measures $(\mu ,\nu )$. In particular, we show that $(\mu ,\nu )$ is a
spectral pair if and only if $(\nu ,\mu )$ is a spectral pair.

After this paper was originally submitted two independent proofs \cite{LRW98}%
, \cite{IoPe98} of Conjecture \ref{Con1.5} have appeared. 

\section{\label{S:construction}Construction of Spectra}

The next two sections are
concerned with the
structure of the discrete sets $\Lambda $
which at the same time
serve as spectra for $I^{d}$
(i.e., the basis property), and
also are sets of vectors $\lambda $
which make the translates $\lambda +I^{d}$
tile $\mathbb{R}^{d}$.

There is a \emph{recursive procedure} for constructing spectral pairs in
higher dimensions from \emph{``factors''} in lower dimension. It is a \emph{%
cross-product} construction, and it applies to any two spectral pairs, $%
\left( \Omega _{i},\Lambda _{i}\right) $, $i=1,2$, in arbitrary dimensions $%
d_{1}$ and $d_{2}$. While it is clear that the ``spectral-pair category'' is 
\emph{closed under tensor product} (see \cite{JoPe92,JoPe94}), the following
result is new:

\begin{theorem}
\label{T:construction}Let $(\Omega _{1},\Lambda _{1})$ be a spectral pair in
dimension $d_{1}$, let $\Omega _{2}$ be a set of finite positive measure in
dimension $d_{2}$. Suppose
that for each 
$\lambda _{1}\in
\Lambda _{1}$,
$\Lambda (\lambda _{1})$
is a
discrete subset of
$\mathbb{R}^{d_{2}}$
such
that $(\Omega _{2},\Lambda (\lambda _{1}))$
is a spectral pair.
If 
$\Lambda =\{(\lambda _{1},\lambda _{2}):\lambda _{1}\in \Lambda _{1},\lambda
_{2}\in \Lambda (\lambda _{1})\}$ then $\left( \Omega _{1}\times \Omega %
_{2},\Lambda \right) $ is a spectral pair in $d_{1}+d_{2}$ dimensions.
\end{theorem}

\begin{proof}%
We first show that the exponentials $\left\{ e_{\lambda }:\lambda \in
\Lambda \right\} $ are mutually orthogonal in $\mathcal{L}%
^{2}\left( \Omega _{1}\times \Omega _{2}\right) $ where the $e_{\lambda }$'s
are given on $\Omega _{1}\times \Omega _{2}$ by the usual formula (\ref{eq1}%
) from Section \ref{S1}. The inner product in $\mathcal{L}^{2}\left( \Omega %
_{1}\times \Omega _{2}\right) $ of $e_{\lambda }$ and $e_{\lambda ^{\prime }}
$ factors as follows: 
\begin{equation*}
\int_{\Omega _{1}}e_{\lambda _{1}-\lambda _{1}^{\prime }}\left( x\right)
\left( \int_{\Omega _{2}}e_{\lambda _{2}-\lambda _{2}^{\prime }}\left(
y\right) \,dy\right) \,dx.
\end{equation*}
If $\lambda _{1}\ne \lambda _{1}^{\prime }$ in $\Lambda _{1}$, then it
vanishes since $\left( \Omega _{1},\Lambda _{1}\right) $ is a spectral pair;
and, if $\lambda _{1}=\lambda _{1}^{\prime }$ but $\lambda _{2}\ne \lambda
_{2}^{\prime }$, it vanishes since $\left( \Omega _{2},\Lambda (\lambda
_{1})\right) $ is one. This proves orthogonality of $\Lambda $. To see that
it is total, let $f\in \mathcal{L}^{2}\left( \Omega _{1}\times \Omega %
_{2}\right) $ and suppose $f$ is orthogonal to $\Lambda $. The inner
products (vanishing) are: 
\begin{equation*}
\ip{e_{\lambda }}{f}%
_{\Omega _{1}\times \Omega _{2}}=\int_{\Omega _{2}}e_{\lambda _{2}}\left(
y\right) e_{\lambda _{1}}\left( y\right) \left( \int_{%
\Omega _{1}}e_{\lambda _{1}}\left( x\right) \overline{f\left( x,y\right) }%
\,dx\right) \,dy.
\end{equation*}
If $\lambda _{1}$ is fixed, and the double integral vanishes for all $%
\lambda _{2}\in \Lambda (\lambda _{1})$, then
the integral
$\int_{\Omega _{1}}e_{\lambda
_{1}}\left( x\right) \overline{f\left( x,y\right) }\,dx=0$ for almost all
$y$, by the totality of $\Lambda (\lambda _{1})$ on $\Omega _{2}$. But
$\lambda_{1}$ is arbitrary so the totality of $\Lambda _{1}$ on $\Omega _{1}$
implies $f=0$. We conclude, that $\Lambda $ is total on $\Omega _{1}\times 
\Omega _{2}$ as claimed.%
\end{proof}%

A more concrete version of Theorem \ref{T:construction} is:

\begin{theorem}
\label{Thm5.3}Let $\left( \Omega _{i},\Lambda _{i}\right) $, $i=1,2$, be
spectral pairs in the respective dimensions $d_{1}$ and $d_{2}$, and let $%
\beta \colon \Lambda _{1}\rightarrow \mathbb{R}^{d_{2}}$ be an arbitrary
function.
Let 
\begin{equation}
\Lambda _{\beta }:=\left\{ 
\begin{pmatrix}
                      \lambda _{1} \\
                      \beta \left( \lambda _{1}\right) +\lambda _{2}
                      \end{pmatrix}
                      %
:\lambda _{1}\in \Lambda _{1}\text{ and }\lambda _{2}\in \Lambda
_{2}\right\} .  \label{eq51}
\end{equation}
Then $\left( \Omega _{1}\times \Omega _{2},\Lambda _{\beta }\right) $ is a
spectral pair in $d_{1}+d_{2}$ dimensions.
\end{theorem}

\begin{proof}%
If $\left( \Omega _{2},\Lambda _{2}\right) $ is a spectral pair, then so is $%
\left( \Omega _{2},\Lambda _{2}+\beta \right) $ for any vector $\beta $. An
application of Theorem \ref{T:construction} completes the proof. 
\end{proof}%

By repeatedly applying Theorem \ref{Thm5.3} if follows that if $\Lambda $ is
the set of points given by: 
\begin{equation}
\begin{pmatrix}
                      \alpha + k_{1} \\
                      \beta_{1} \left( k_{1}\right) +k_{2} \\
                      \beta_{2} \left( k_{1},k_{2}\right) +k_{3} \\
                      \vdots  \\
                      \beta_{d-1} \left( k_{1},k_{2},\dots ,k_{d-1}\right) +k_{d}
                      \end{pmatrix}
                      %
\label{eq17}
\end{equation}
with $k_{1},k_{2},\dots ,k_{d}\in \mathbb{Z}$, where 
$\beta_{i}  \colon \mathbb{Z}^{i}\longrightarrow \left[ 0,1\right\rangle $
are fixed functions, then $(I^{d},\Lambda )$ is a spectral pair. Of course,
there are the obvious modifications resulting from permutation of the $d$
coordinates; but, when $d\ge 10$, these configurations do \emph{not} suffice
for cataloguing all the possible spectra $\Lambda $ which turn $\left(
I^{d},\Lambda \right) $ into an $\mathbb{R}^{d}$-spectral pair;
see Section \ref{S:periodic}.

We now turn to a result which is a
partial
converse to Theorem \ref{T:construction},
its statement requires the following notation. It is motivated by the
``projection method'' from quasicrystal theory; see, e.g., \cite{Hof95}. For
subsets $\Lambda $ of the Cartesian product $\mathbb{R}^{d_{1}}\times %
\mathbb{R}^{d_{2}}$, let 
\begin{equation*}
\Lambda _{1}:=P_{1}\Lambda =\left\{ \lambda _{1}\in \mathbb{R}%
^{d_{1}}:\exists \lambda _{2}\in \mathbb{R}^{d_{2}}\text{ s.t.\ }%
\begin{pmatrix}
       \lambda _{1} \\
       \lambda _{2} 
       \end{pmatrix}%
\in \Lambda \right\} 
\end{equation*}
(where it is convenient here to use column vector formalism). We shall also
need sections of $\Lambda $ in the $\mathbb{R}^{d_{2}}$ coordinate direction
as follows: If $\lambda _{1}\in \Lambda _{1}$ ($=P_{1}\Lambda $), set 
\begin{equation*}
\Lambda \left( \lambda _{1}\right) :=\left\{ \lambda _{2}\in \mathbb{R}%
^{d_{2}}:%
\begin{pmatrix}
       \lambda _{1} \\
       \lambda _{2} 
       \end{pmatrix}%
\in \Lambda \right\} .
\end{equation*}

\begin{lemma}
\label{Lem5.4}Let $\Omega _{i}\subset \mathbb{R}^{d_{i}}$, $i=1,2$, be
subsets with finite positive Lebesgue measure in the respective dimensions,
and let $\Lambda \subset \mathbb{R}^{d_{1}}\times \mathbb{R}^{d_{2}}$ be a
subset such that $\left( \Omega _{1}\times \Omega _{2},\Lambda \right) $ is
a spectral pair in $d_{1}+d_{2}$ dimensions. Then for every $\lambda _{1}\in
P_{1}\Lambda $, the exponentials $\left\{ e_{\xi }^{\left( 2\right) }:\xi
\in \Lambda \left( \lambda _{1}\right) \right\} $ are orthogonal in $%
\mathcal{L}^{2}\left( \Omega _{2}\right) $; and they are total in $\mathcal{L%
}^{2}\left( \Omega _{2}\right) $ if and only if $e_{\lambda _{1}}^{\left(
1\right) }$ and $e_{\lambda _{1}^{\prime }}^{\left( 1\right) }$ are
orthogonal in $\mathcal{L}^{2}\left( \Omega _{1}\right) $ for all $\lambda
_{1}^{\prime }\in P_{1}\left( \Lambda \right) \diagdown \left\{ \lambda
_{1}\right\} $.
\end{lemma}

\begin{proof}%
To check orthogonality, let $\xi ,\eta \in \Lambda \left( \lambda
_{1}\right) $. Then the two points $\left( 
\begin{smallmatrix}
\lambda _{1} \\
\vphantom{\eta}\smash{\xi} 
\end{smallmatrix}\right) $ and $\left( 
\begin{smallmatrix}
\lambda _{1} \\
\eta 
\end{smallmatrix}\right) $ are in $\Lambda $, and the corresponding $\Omega
_{1}\times \Omega _{2}$-inner product is zero. But it is also 
\begin{equation*}
m_{d_{1}}\left( \Omega _{1}\right) 
\ip{e_{\xi }}{e_{\eta }}%
_{\Omega _{2}},
\end{equation*}
and since $m_{d_{1}}\left( \Omega _{1}\right) >0$, the orthogonality follows.

We now show that $\Lambda \left( \lambda _{1}\right) $ is total in $\mathcal{%
L}^{2}\left( \Omega _{2}\right) $ if $%
\ip{e^{\left( 1\right) }_{\lambda _{1}}}{e^{\left( 1\right) }_{\lambda ^{\prime }_{1}}}%
_{\Omega _{1}}=0$ for all $\lambda _{1}^{\prime }\in P_{1}\Lambda $,
$\lambda _{1}^{\prime }\neq \lambda _{1}^{{}}$.
Let $%
g\in \mathcal{L}^{2}\left( \Omega _{2}\right) $ and suppose $g$ is
orthogonal to all the $\Lambda \left( \lambda _{1}\right) $-exponentials.
Let $\left( 
\begin{smallmatrix}
\lambda ^{\prime }_{1} \\
\lambda ^{\prime }_{2} 
\end{smallmatrix}\right) $ be a general point in $\Lambda $. Then the inner
product with $e_{\lambda _{1}}^{\left( 1\right) }\otimes g$ is 
\begin{equation*}
\ip{e^{\left( 1\right) }_{\lambda ^{\prime }_{1}}}{e^{\left( 1\right) }_{\lambda _{1}}}%
_{\Omega _{1}} 
\ip{e^{\left( 2\right) }_{\lambda ^{\prime }_{2}}}{g}%
_{\Omega _{2}}.
\end{equation*}
If $\lambda _{1}^{\prime }=\lambda _{1}$, then $\lambda _{2}^{\prime }\in
\Lambda \left( \lambda _{1}\right) $, and the second factor vanishes. If $%
\lambda _{1}^{\prime }\ne \lambda _{1}$, then the first factor is zero, and
we get that $e_{\lambda _{1}}^{\left( 1\right) }\otimes g$ is orthogonal in $%
\mathcal{L}^{2}\left( \Omega _{1}\times \Omega _{2}\right) $ to the $\Lambda 
$-exponentials. They are total, and we conclude that $g$ vanishes in $%
\mathcal{L}^{2}\left( \Omega _{2}\right) $.

The remaining case is when $%
\ip{e^{\left( 1\right) }_{\lambda _{1}}}{e^{\left( 1\right) }_{\lambda ^{\prime }_{1}}}%
_{\Omega _{1}}\ne 0$ for some $\lambda _{1}^{\prime }\in P_{1}\left( \Lambda
\right) \diagdown \left\{ \lambda _{1}\right\} $. But it follows that then $%
\ip{e^{\left( 2\right) }_{\xi }}{e^{\left( 2\right) }_{\eta }}%
_{\Omega _{2}}=0$ for all $\xi \in \Lambda \left( \lambda _{1}\right) $ and $%
\eta \in \Lambda \left( \lambda _{1}^{\prime }\right) $. In particular, $%
\Lambda \left( \lambda _{1}\right) $ is then \emph{not} total in $\mathcal{L}%
^{2}\left( \Omega _{2}\right) $.%
\end{proof}%

\section{Dimensions Two and Three\label{S2}}

In this section we prove Conjecture \ref{Con1.5} for $d=1,2,3$. Furthermore
we give a complete classification of the possible spectra for the
unit cube in those dimensions.

We begin with the following simple observation in one dimension for $\Omega
=I=\left[ 0,1\right\rangle $.

\begin{proposition}
\label{Pro2.1}The only subsets $\Lambda \subset \mathbb{R}$ such that $%
\left( I,\Lambda \right) $ is a spectral pair are the translates 
\begin{equation}
\Lambda _{\alpha }:=\alpha +\mathbb{Z}=\left\{ \alpha +n:n\in \mathbb{Z}%
\right\}   \label{eq13}
\end{equation}
where $\alpha $ is some fixed real number.
\end{proposition}

In two dimensions, the corresponding result is more subtle, but the
possibilities may still be enumerated as follows:

\begin{theorem}
\label{Thm2.2}The only subsets $\Lambda \subset \mathbb{R}^{2}$ such that $%
\left( I^{2},\Lambda \right) $ is a spectral pair must belong to either one
or the other of the two classes, indexed by a number $\alpha $, and a
sequence $\left\{ \beta _{m}\in \left[ 0,1\right\rangle :m\in \mathbb{Z}%
\right\} $, where 
\begin{align}
\Lambda & =\left\{ 
\begin{pmatrix}
                      \alpha +m \\
                      \beta _{m}+n
                      \end{pmatrix}
                      %
:m,n\in \mathbb{Z}\right\}   \label{eq14} \\
\intertext{or}%
\Lambda & =\left\{ 
\begin{pmatrix}
                      \beta _{n}+m \\
                      \alpha +n
                      \end{pmatrix}
                      %
:m,n\in \mathbb{Z}\right\} .  \label{eq15}
\end{align}
Each of the two types occurs as the spectrum of a pair for the cube $I^{2}$,
and each of the sets $\Lambda $ as specified is a tiling set for the cube $%
I^{2}$.
\end{theorem}

\begin{proof}%
The assertion in the theorem about $\Lambda $-translations tiling the plane
with $I^{2}$ is clear from (\ref{eq14})--(\ref{eq15}), and it is illustrated
graphically in Figure \ref{tiling1}. 

\begin{figure}[tbp]
\makebox[\textwidth]{%
\psfig{file=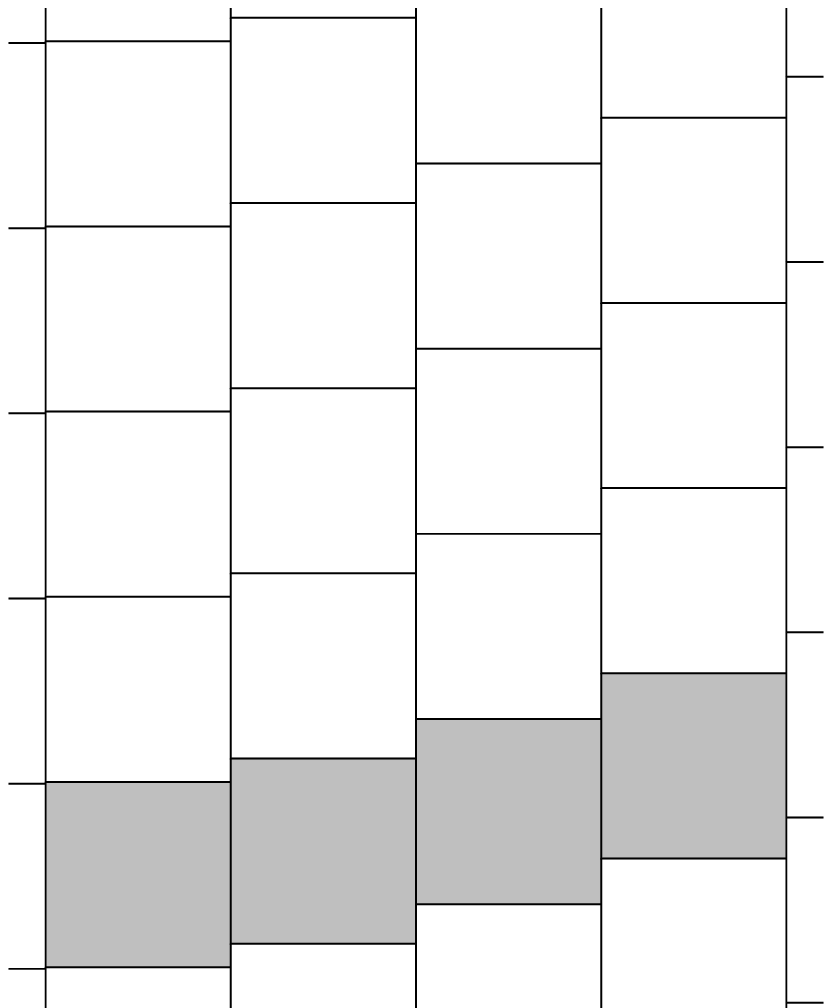,bbllx=26bp,bblly=0bp,bburx=262bp,bbury=288bp,height=187pt}\hfill%
\raisebox{18pt}%
{\psfig{file=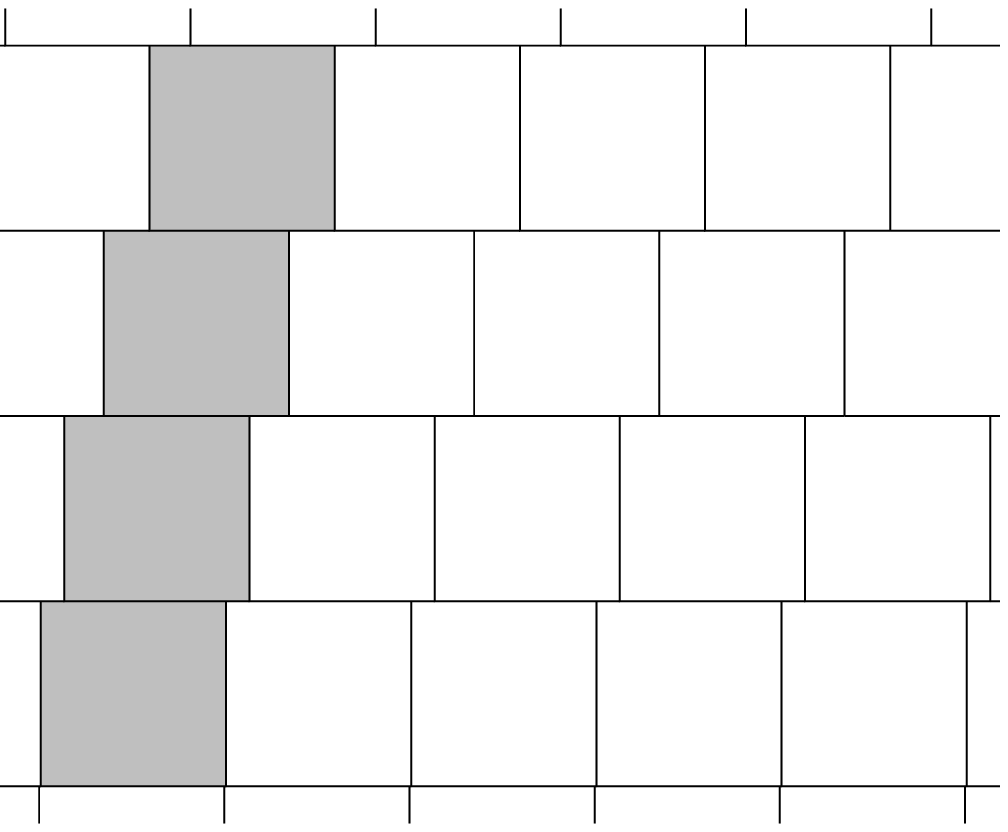,bbllx=0bp,bblly=26bp,bburx=288bp,bbury=262bp,width=187pt}}}%
\caption{Illustrating tiling with (\ref{eq14}) (left) and (\ref{eq15}) (right)}
\label{tiling1}
\end{figure}
%

It is immediate from Theorem \ref{Thm5.3} that each one of the two formulas (%
\ref{eq14})--(\ref{eq15}) for $\Lambda $ make $\left( I^{2},\Lambda \right) $
a spectral pair, and the main result is that there are not others. We show
this directly by an examination of the possibilities for $\Lambda $ which
are implied by the inclusion 
\begin{equation}
\Lambda -\Lambda \subset \mathbf{Z}_{I^{2}}
\cup \left\{ 0\right\}  \label{eq21}
\end{equation}
where $\mathbf{Z}_{I^{2}}$ is read off from Lemma \ref{Lem1.4} above. Again
a translation of $\Lambda $ by a single vector in the plane will reduce the
analysis to the case when $\left( 
\begin{smallmatrix}
0 \\
0
\end{smallmatrix}\right) $ is in $\Lambda $.
Let $\lambda =\left( 
\xi  ,
\eta 
\right) \in \Lambda $
and suppose $\lambda \notin \mathbb{Z}^{2}$.
Then either $\xi $ or $\eta $ is not an integer.
Suppose $\eta $ is not an integer.
Then $\xi $ is a nonzero integer.
Let $\lambda^{\prime}=\left( 
\xi^{\prime}  ,
\eta^{\prime}
\right) $
be an arbitrary point in $\Lambda $.
If $\xi
^{\prime }$ is not an integer, then $\eta^{\prime}$
is a nonzero integer, so $\lambda -\lambda^{\prime}$
is not in $\mathbf{Z}_{I^{2}}$, contradicting (\ref{eq21}).
So $\xi^{\prime }\in\mathbb{Z}$
for any $\lambda^{\prime}=\left( 
\xi^{\prime}  ,
\eta^{\prime}
\right) $ in $\Lambda $.
To verify $\Lambda $ is a subset
of a set given by (\ref{eq14}) for $\alpha =0$, we need only check that if $%
\left( \begin{smallmatrix}
\xi  \\
\eta 
\end{smallmatrix}\right) $ and $\left( 
\begin{smallmatrix}
\xi ^{\prime} \\
\eta ^{\prime}
\end{smallmatrix}\right) $ are different points in $\Lambda $, with $\xi
=\xi ^{\prime }\in \mathbb{Z}$, then $\eta -\eta ^{\prime }\in \mathbb{Z}$.
But recall $\left( 
\begin{smallmatrix}
0 \\
\eta -\eta ^{\prime }
\end{smallmatrix}\right) \in \mathbf{Z}_{I^{2}}$, so it follows from Lemma 
\ref{Lem1.4} that $\eta -\eta ^{\prime }\in \mathbb{Z}\diagdown \left\{
0\right\} $. Since an orthonormal basis cannot be a strict subset of another
orthonormal basis for the same space it follows that $\Lambda $ is given by (%
\ref{eq14}).%
\end{proof}%

Replacing the appeal to Lemma \ref{Lem1.4} in this proof with an appeal to
Theorem \ref{T:Keller} it follows that any tiling set $\Lambda $ for the
cube $I^{2}$ must be given by (\ref{eq14})--(\ref{eq15}), we leave the
details for the reader. The fact that this simple tiling pattern for the
cube $I^{d}$ in $d$ dimensions is broken for $d=10$ follows from examples of
Lagarias and Shor \cite{LaSh92}. It is shown there that for each $d\ge 10$
there exists a tiling of $\mathbb{R}^{d}$ by translates of $I^{d}$ such that
no two tiles have a complete facet in common. These examples also
demonstrate, see Section \ref{S:periodic}, that if $d\ge 10$, then the
corresponding combinations (\ref{eq17})
do not supply all possible
spectra
for $I^{d}$.

The following result
shows that spectra for $I^{3}$ and tilings of $\mathbb{R}^{3}$ by $I^{3}$
are the same by fully determining each.
No complete description of such tilings or spectra is known for
$d>3$.

\begin{theorem}
\label{Thm5.5}$\left( I^{3},\Lambda \right) $ 
is a tiling pair, or a spectral pair, if and only if, 
after a possible translation by a single vector 
and a possible permutation of the coordinates $\left(
x_{1},x_{2},x_{3}\right) $, $\Lambda $ can be brought into the following
form: there is a partition of $\mathbb{Z}$ into disjoint subsets $A$, $B$
\textup{(}one possibly empty\/\textup{)} with
associated functions 
\begin{align*}
\alpha _{0}\colon A& \longrightarrow \left[ 0,1\right\rangle , \\
\alpha _{1}\colon A\times \mathbb{Z}& \longrightarrow \left[
0,1\right\rangle , \\
\beta _{0}\colon B& \longrightarrow \left[ 0,1\right\rangle , \\
\beta _{1}\colon B\times \mathbb{Z}& \longrightarrow \left[ 0,1\right\rangle 
\end{align*}
such that $\Lambda $ is the \textup{(}disjoint\textup{)} union of 
\begin{equation}
\begin{pmatrix}
                      a \\
                      \alpha _{0}\left( a\right) +k \\
                      \alpha _{1}\left( a,k\right) +l
                      \end{pmatrix}
                      %
\text{\quad and\quad }%
\begin{pmatrix}
                      b \\
                      \beta _{1}\left( b,n\right) +m \\
                      \beta _{0}\left( b\right) +n
                      \end{pmatrix}
                      %
\label{E:3d}
\end{equation}
as $a\in A$, $b\in B$, and $k,l,m,n\in \mathbb{Z}$.
\end{theorem}

\begin{proof}%
Suppose
$\Lambda$ is a tiling set for
$I^3$. By the tiling property there
exist
functions $\alpha ,
\beta , \gamma \colon \mathbb{Z}^3\mapsto\left\langle-1,0\right]$ so that 
\begin{equation*}
   \Lambda=\left\{
      \begin{pmatrix}
         n+\alpha(l,m,n) \\
         m+\beta(l,m,n) \\
         l+\gamma(l,m,n)
      \end{pmatrix}
      : l,m,n\in\mathbb{Z}\right\} .
\end{equation*}
Fix $l,m\in\mathbb{Z}$ then Theorem
\ref{T:Keller}
(Keller's theorem) implies
$\alpha(l,m,n)$ is independent of $n$, we will write
$\alpha(l,m)$ in place of $\alpha(l,m,n)$ to indicate this 
independence. Similarly, $\beta(l,m,n)=\beta(l,n)$ and
$\gamma(l,m,n)=\gamma(m,n)$. 
Considering, for fixed $n\in\mathbb{Z}$, 
the intersection of the plane
$x_1=n$ by cubes $I^3+\lambda$, $\lambda\in\Lambda$ it follows that 
the set
\begin{equation*}
   \Lambda_{1,n}=\left\{
      \begin{pmatrix}
         m+\beta(l,n) \\
         l+\gamma(m,n)
      \end{pmatrix}
      : l,m\in\mathbb{Z}\right\}
\end{equation*}
is a tiling set for $I^2$ in $\mathbb{R}^2$. Hence, by our $d=2$
result (Theorem
\ref{Thm2.2}),
either
\begin{equation*}
   \Lambda_{1,n}=\left\{
      \begin{pmatrix}
         m+\tilde{\beta}(l,n) \\
         l+\tilde{\gamma}(n)
      \end{pmatrix}
      : l,m\in\mathbb{Z}\right\}
   \text{\qquad or\qquad }
   \Lambda_{1,n}=\left\{
      \begin{pmatrix}
         m+\tilde{\beta}(n) \\
         l+\tilde{\gamma}(m,n)
      \end{pmatrix}
      : l,m\in\mathbb{Z}\right\}.\notag
\end{equation*}
It follows that there exist $A,B\subset\mathbb{Z}$ so that 
$A\cup B=\mathbb{Z}$, $A\cap B=
\varnothing$,
and 
\begin{equation*}
   \Lambda=
      \left\{
      \begin{pmatrix}
         n+\alpha(l,m) \\
         m+\beta(l,n) \\
         l+\gamma(n)
      \end{pmatrix}
      : l,m\in\mathbb{Z},n\in A\right\}
   \cup
      \left\{
      \begin{pmatrix}
         n+\alpha(l,m) \\
         m+\beta(n) \\
         l+\gamma(m,n)
      \end{pmatrix}
      : l,m\in\mathbb{Z},n\in B\right\}.\notag
\end{equation*}
For each $m\in\mathbb{Z}$ 
\begin{equation*}
   \Lambda_{2,m}=\left\{
      \begin{pmatrix}
         n+\alpha(l,m) \\
         l+\gamma(n)
      \end{pmatrix}
      : l\in\mathbb{Z}, n\in A\right\}
   \cup
      \left\{
      \begin{pmatrix}
         n+\alpha(l,m) \\
         l+\gamma(m,n)
      \end{pmatrix}
      : l\in\mathbb{Z}, n\in B\right\}
\end{equation*}
is a tiling set for $I^2$ in $\mathbb{R}^2$. Hence, our $d=2$
result implies either (1) $\alpha(l,m)=\alpha(m)$ for all
$l,m\in\mathbb{Z}$ or (2) $\gamma(m,n)=\gamma(n)$, for all
$m\in\mathbb{Z}$, $n\in B$. Suppose (1), then for each
$l\in\mathbb{Z}$ 
\begin{equation*}
   \Lambda_{3,l}=\left\{
      \begin{pmatrix}
         n+\alpha(m) \\
         m+\beta(l,n)
      \end{pmatrix}
      : m\in\mathbb{Z}, n\in A\right\}
   \cup
      \left\{
      \begin{pmatrix}
         n+\alpha(m) \\
         m+\beta(n)
      \end{pmatrix}
      : m\in\mathbb{Z}, n\in B\right\}
\end{equation*}
is a tiling set for $I^2$, hence our $d=2$ result implies that
either (1a) $\alpha(m)=\alpha_0$ for all $m\in\mathbb{Z}$ or (1b)
$\beta(l,n)=\beta(l)$ for $n\in A$ and $\beta(n)=\beta_0$ for 
$n\in B$. If (1a) then we are done, so suppose (1b): then
\begin{equation*}
   \Lambda=
      \left\{
      \begin{pmatrix}
         n+\alpha(m) \\
         m+\beta(l) \\
         l+\gamma(n)
      \end{pmatrix}
      : l,m\in\mathbb{Z},n\in A\right\}
   \cup
      \left\{
      \begin{pmatrix}
         n+\alpha(m) \\
         m+\beta_0 \\
         l+\gamma(m,n)
      \end{pmatrix}
      : l,m\in\mathbb{Z},n\in B\right\}.\notag
\end{equation*}
Let, if possible, $n_1\in A$, $n_2\in B$, $m_1,m_2,l\in\mathbb{Z}$
be
such
that $\alpha(m_1)\neq\alpha(m_2)$ and
$\beta(l)\neq\beta_0$: then 
\begin{equation*}
      \begin{pmatrix}
         n_1+\alpha(m_1) \\
         m_1+\beta(l) \\
         l+\gamma(n_1)
      \end{pmatrix}
      -
       \begin{pmatrix}
         n_2+\alpha(m_2) \\
         m_2+\beta_0 \\
         l+\gamma(n_2)
      \end{pmatrix}
\end{equation*}
does not have any nonzero integer entry, contradicting Keller's
theorem. So, either $B=
\varnothing$,
$A=
\varnothing$,
$\alpha(m)=\alpha_0$ for all $m\in\mathbb{Z}$, or $\beta(l)=\beta_0$
for all $l\in\mathbb{Z}$; in the three last cases we are done, so
assume $A=\mathbb{Z}$. If $\alpha(m_1)\neq\alpha(m_2)$,
$\beta(l_1)\neq\beta(l_2)$, and $\gamma(n_1)\neq\gamma(n_2)$ then 
\begin{equation*}
      \begin{pmatrix}
         n_1+\alpha(m_1) \\
         m_1+\beta(l_1) \\
         l_1+\gamma(n_1)
      \end{pmatrix}
      -
       \begin{pmatrix}
         n_2+\alpha(m_2) \\
         m_2+\beta(l_2) \\
         l_2+\gamma(n_2)
      \end{pmatrix} 
\end{equation*} 
does not have any nonzero integer
entry, contradicting Keller's theorem. This contradiction completes the
proof of case (1). The proof of case (2) is similar; we leave the details
for the reader.

Conversely, for every such set $\Lambda $, the translates of $I^{3}$ by the
vectors of $\Lambda $ clearly tile $\mathbb{R}^{3}$. This completes the
description of tilings of $\mathbb{R}^{3}$ by $I^{3}$.

By Theorem \ref{T:construction} and Theorem \ref{Thm2.2} any set $\Lambda $
of the form (\ref{E:3d}) is a spectrum for $I^{3}$.
We sketch a proof of the converse.

Apply Lemma \ref{Lem5.4} to $\Omega _{1}%
\times \Omega _{2}=I\times I^{2}$,
to
show that if $\left( I^{3},\Lambda \right) $ is a spectral pair, then one of
the three coordinate intervals may be picked as $\Omega _{1}$
in Lemma \ref{Lem5.4}, i.e., $\Omega %
_{1}=I$, $\Omega _{2}=I^{2}=I\times I$ and with $P_{1}\Lambda =\Lambda _{1}$
satisfying the orthogonality on $\mathcal{L}^{2}\left( I\right) $. By
Lemma \ref{Lem1.4}, 
this means 
\begin{equation}
\Lambda _{1}-\Lambda _{1}\subset \mathbb{Z}.  \label{eq53}
\end{equation}
Eventually we show that $\Lambda _{1}$ must be of the form $\theta _{1}+%
\mathbb{Z}$. But to select the one of the three coordinates which has this
form, consider the canonical mapping 
\begin{equation*}
\mathbb{R}\overset{\pi }{\longrightarrow }\mathbb{R}\diagup \mathbb{Z}\simeq
\left[ 0,1\right\rangle 
\end{equation*}
and select the one of the three sets $\pi \left( P_{j}\left( \Lambda \right)
\right) $, $j=1,2,3$, of the smallest cardinality. Assume it is $j=1$ for
simplicity. We will show that $\Lambda _{1}$ then satisfies (\ref{eq53}), so
that Lemma \ref{Lem5.4} applies. The assertion is that the set $\pi \left(
\Lambda _{1}\right) $ is a singleton. The proof is indirect. Suppose \emph{%
ad absurdum} that, for some $\theta _{1}$ such that $0<\theta _{1}<1$, $%
\Lambda _{1}$ meets both $\mathbb{Z}$ and $\theta _{1}+\mathbb{Z}$. We
conclude
from Theorem \ref{Thm2.2}
that for $\lambda _{1}$ in each of the two sets $\mathbb{Z}$ or $%
\theta _{1}+\mathbb{Z}$, the points in $\Lambda \left( \lambda _{1}\right) $
must be of the form $\left( 
\begin{smallmatrix}
k \\
\alpha \left( k\right) +l 
\end{smallmatrix}\right) $ for $k,l\in \mathbb{Z}$ and $\alpha \colon %
\mathbb{Z}\rightarrow \left[ 0,1\right\rangle $ some function, or
alternatively $\left( 
\begin{smallmatrix}
\beta \left( l^{\prime }\right) +k^{\prime } \\
l^{\prime }
\end{smallmatrix}\right) $, $k^{\prime },l^{\prime }\in \mathbb{Z}$ with $%
\beta $ some possibly different function.
Calculating $\mathcal{L}^{2}\left(
I^{3}\right) $-inner products for associated points $\lambda ,\lambda
^{\prime }\in \Lambda $ with $P_{1}\left( \lambda \right) =m\in \mathbb{Z}$,
and $P_{1}\left( \lambda ^{\prime }\right) =n+\theta _{1}$ ($n\in \mathbb{Z}$%
), we get the following possibilities for the respective coordinates in the
second and third place: 
\begin{equation*}
\begin{pmatrix}
                      \alpha \left( m\right) +k \\
                      \beta \left( m,k\right) +l
                      \end{pmatrix}
                      %
\text{\quad or\quad }%
\begin{pmatrix}
                      \delta \left( m,l^{\prime }\right) +k^{\prime } \\
                      \gamma \left( m\right) +l^{\prime }
                      \end{pmatrix}
                      %
.
\end{equation*}
But if the difference $P_{1}\lambda -P_{1}\lambda ^{\prime }$ is not in $%
\mathbb{Z}$, then one of the corresponding differences in the second place,
or the third place, must be in $\mathbb{Z}$. Making variations, we conclude
that then one of the two sets $\pi \left( P_{j}\Lambda \right) $, $j=2,3$,
must be a singleton. But this contradicts that $\pi \left( P_{1}\Lambda
\right) $ has \emph{two} distinct points, and is chosen to have smallest
cardinality of the three sets $\pi \left( P_{j}\Lambda \right) $, $j=1,2,3$.%
\end{proof}%

\begin{corollary}
\label{CorNew3.4}The commuting
self-adjoint extensions
$\left\{ H_{j}:j=1,\dots ,d\right\}$ in \textup{(\ref{eq7})} are
completely classified and determined, for
$d=1,2,3$ and $\Omega =I^{d}$, by
Proposition \textup{\ref{Pro2.1}} for $d=1$,
Theorem \textup{\ref{Thm2.2}} for $d=2$, and
Theorem \textup{\ref{Thm5.5}} for $d=3$.
\end{corollary}

\begin{proof}
The stated conclusion
follows from combining the results
in the present section with
Theorem \ref{Thm1.1}; for $I^{d}$ the spectral
condition is equivalent to the operator
extension property.
\end{proof}

\section{\label{S:periodic}Periodic Sets}

A discrete set $T\subset \mathbb{R}^{d}$ is \emph{periodic}
if there exists a finite set $L\subset \mathbb{R}^{d}$
and an invertible $d\times d$ matrix $R$
with real entries such that $T=L+R\mathbb{Z}^{d}$.
Periodic sets have played an important
role in the study of spectral pairs,
see, e.g., 
\cite{JoPe92}, \cite{Ped96}, \cite{LaWa97b},
and in the study of tilings
by translation, see, e.g., \cite{LaSh92}, \cite{LaWa96c}.
In this section we establish Conjecture \ref{Con1.5}
under the further hypothesis that $L$ is a
periodic set.
The periodic case is of
interest because (\cite{Kel30}) Keller's conjecture
about cube tilings (see below) is false if and
only if it is false for certain periodic tilings.
Also a long-standing conjecture is that any
bounded tile in $\mathbb{R}^d$ admits a
periodic tiling set.
 
\begin{lemma}
\label{Lem6.1}Let $R$ be an invertible $d\times d$ matrix with real entries
and let $\Lambda :=R\mathbb{Z}^{d}$ be the corresponding lattice. If $\Omega
\subset \mathbb{R}^{d}$ is a measurable set, then $\Omega $ is a $\Lambda $-tile
if and only if

\begin{enumerate}
\item  \label{Lem6.1(1)}$m\left( \Omega \right) =\lvert \det R\rvert$ and

\item  \label{Lem6.1(2)}$m\left( \Omega \cap \left( \Omega +l\right) \right)
=0$ for all $l\in \Lambda \diagdown \left\{ 0\right\} $,
\end{enumerate}

\noindent where $m$ is the $d$-dimensional Lebesgue measure.
\end{lemma}

A
measurable set $\Omega $ is called a $\Lambda $\emph{-tile} if $%
\bigcup_{l\in \Lambda }\left( \Omega +l\right) $ is a measure-theoretic partition
of $\mathbb{R}^{d}$. Note that $\Omega _{R}:=RI^{d}$ is a $\Lambda $-tile.

\begin{proof}[Proof of Lemma \textup{\ref{Lem6.1}}]
If $\Omega $ is a $\Lambda $-tile then $m\left( \Omega \right) =m\left( \Omega
_{R}\right) =\lvert \det R\rvert$, since any two $\Lambda $-tiles have the same
volume \cite{GrLe87}. Conversely, suppose a measurable set $\Omega $ has
properties (\ref{Lem6.1(1)}) and (\ref{Lem6.1(2)}). Let $\Omega _{l}:=\left(
\Omega _{R}+l\right) \cap \Omega $; then $\bigcup_{l\in \Lambda }\Omega _{l}$ is a
measure-theoretic partition of $\Omega $. By (\ref{Lem6.1(2)}), the sets $%
\Omega _{l}-l=\Omega _{R}\cap \left( \Omega -l\right) $, $l\in \Lambda $, are
measure disjoint. Hence, 
\begin{equation*}
m\left( 
\smash{\bigcup _{l\in \Lambda }}\vphantom{\bigcup }
\left( \Omega _{l}-l\right) \right) =m\left( 
\smash{\bigcup _{l\in \Lambda }}\vphantom{\bigcup }
\left( \Omega _{l}\right) \right) =m\left( \Omega \right) =\lvert \det R%
\rvert
\vphantom{\bigcup _{l\in \Lambda }}
\end{equation*}
by (\ref{Lem6.1(1)}). It follows that $\bigcup_{l\in \Lambda }\left( \Omega
_{l}-l\right) =\Omega _{R}$, up to sets of measure zero. Hence, up to sets
of measure zero, 
\begin{multline*}
\bigcup_{k\in \Lambda }\left( \Omega +k\right)  =\bigcup_{k\in \Lambda }\bigcup_{l\in
\Lambda }\left( \Omega _{l}+k\right) 
 =\bigcup_{l\in \Lambda }\bigcup_{k\in \Lambda }\left( \Omega _{l}+k\right) \\
 =\bigcup_{l\in \Lambda }\bigcup_{k\in \Lambda }\left( \Omega _{l}-l+k\right) 
 =\bigcup_{k\in \Lambda }\bigcup_{l\in \Lambda }\left( \left( \Omega _{l}-l\right)
+k\right) 
 =\bigcup_{k\in \Lambda }\left( \Omega _{R}+k\right) 
 =\mathbb{R}^{d},
\end{multline*}
as we needed to show.%
\end{proof}%

\begin{theorem}
\label{Thm6.5}Let $R$ be an invertible $d\times d$ matrix with real entries,
let $\mathbf{L}\subset \mathbb{R}^{d}$ be finite
with $\left| \mathbf{L}\right| $ elements,
and let $L:=R\mathbb{Z}%
^{d} $ be the lattice in $\mathbb{R}^{d}$ generated by the columns of $R$.
If $L\cap \left( \mathbf{L}-\mathbf{L}\right) =\left\{ 0\right\} $, then the
following are equivalent:

\begin{enumerate}
\item  \label{Thm6.5(1)}$\left( I^{d},L+\mathbf{L}\right) $ is a spectral
pair.

\item  \label{Thm6.5(2)}$I^{d}$ is an $\left( L+\mathbf{L}\right) $-tile.

\item  \label{Thm6.5(3)}$\lvert \mathbf{L}\rvert=\lvert \det R\rvert$ and $%
L+\left( \mathbf{L}-\mathbf{L}\right) \subset \mathbf{Z}_{I^{d}}\cup \left\{
0\right\} $.
\end{enumerate}
\end{theorem}

\begin{proof}%
(\ref{Thm6.5(1)}) $\Rightarrow $ (\ref{Thm6.5(3)}) Since the functions $%
\left( e_{\lambda }\right) _{\lambda \in L+\mathbf{L}}$ are orthogonal it
follows that $L+\left( \mathbf{L}-\mathbf{L}\right) \subset \mathbf{Z}%
_{I^{d}}\cup \left\{ 0\right\} $. By \cite[Theorem 4.2]{Lan67} $L+\mathbf{L}$
has uniform density $1$, hence $\lvert \mathbf{L}\rvert=\lvert \det R\rvert$%
. (\ref{Thm6.5(3)}) $\Rightarrow $ (\ref{Thm6.5(1)}) Using $L+\left( \mathbf{%
L}-\mathbf{L}\right) \subset \mathbf{Z}_{I^{d}}\cup \left\{ 0\right\} $ we
conclude that the functions $\left( e_{\lambda }\right) _{\lambda \in L+%
\mathbf{L}}$ are orthogonal in $\mathcal{L}^{2}\left( I^{d}\right) $. So,
since $\lvert \mathbf{L}\rvert=\lvert \det R\rvert$, an application of 
\cite[Theorem 1]{Ped96} allows us to conclude $\left( I^{d},L+\mathbf{L}%
\right) $ is a spectral pair. (\ref{Thm6.5(2)}) $\Rightarrow $ (\ref
{Thm6.5(3)}) That $L+\left( \mathbf{L}-\mathbf{L}\right) \subset \mathbf{Z}%
_{I^{d}}\cup \left\{ 0\right\} $ follows from $I^{d}$ being an $L+\mathbf{L}$%
-tile, by Theorem \ref{T:Keller}. If $l\ne l^{\prime }$ in $\mathbf{L}$,
then $l-l^{\prime }$ has a nonzero integer entry, hence $m\left( \left(
I^{d}+l\right) \cap \left( I^{d}+l^{\prime }\right) \right) =0$. Hence $%
I^{d}+\mathbf{L}$ has volume $\lvert \mathbf{L}\rvert$. By assumption $I^{d}+%
\mathbf{L}$ is an $L$-tile, so using Lemma \ref{Lem6.1} we conclude $\lvert 
\mathbf{L}\rvert=\lvert \det R\rvert$. (\ref{Thm6.5(3)}) $\Rightarrow $ (\ref
{Thm6.5(2)}) Let $x,x^{\prime }\in L$, $l,l^{\prime }\in \mathbf{L}$. If $%
x\ne x^{\prime }$ or $l=l^{\prime }$, then $\left( l-l^{\prime }\right)
+\left( x-x^{\prime }\right) \in \mathbf{Z}_{I^{d}}$ so 
\begin{equation}
m\left( \left( I^{d}+x+l\right) \cap \left( I^{d}+x^{\prime }+l^{\prime
}\right) \right) =0.  \label{eq6.1}
\end{equation}
In particular $I^{d}+\mathbf{L}$ has volume $\lvert \mathbf{L}\rvert=\lvert %
\det R\rvert$ and $m\left( \left( I^{d}+\mathbf{L}\right) \cap \left( I^{d}+%
\mathbf{L}+x\right) \right) =0$ for all $x\in L\diagdown \left\{ 0\right\} $%
. Hence $I^{d}+\mathbf{L}$ is an $L$-tile by Lemma \ref{Lem6.1}. Using (\ref
{eq6.1}) it follows that $I^{d}$ is an $\left( L+\mathbf{L}\right) $-tile.%
\end{proof}%

While studying a conjecture of Minkowski, Keller \cite{Kel30,Kel37} made the
stron\-ger conjecture that, if $T\subset \mathbb{R}^{d}$ is any set such that $%
I^{d}$ is a $T$-tile, then $T-T$ contains one of the canonical unit basis
vectors for $\mathbb{R}^{d}$. Keller's conjecture is true for $d\le 6$ \cite
{Per40} (see also \cite{StSz94}), and false for $d\ge 10$ \cite{LaSh92}. It
remains open for $d=7,8,9$. The $d=10$ example in \cite{LaSh92} has a set of
translations of the form $T=\left( 2\mathbb{Z}\right) ^{10}+\mathbf{L}$,
where $\mathbf{L}$ has $2^{10}=1024$ elements. To check (\ref{Thm6.5(3)}) in
Theorem \ref{Thm6.5} it is sufficient to check that 
\begin{equation*}
\mathbf{L}-\mathbf{L}\subset \mathbf{Z}_{I^{10}}\cup \left\{ 0\right\} ,
\end{equation*}
but this is condition (a) in \cite[p.\ 280]{LaSh92}, hence it is satisfied
by construction. It follows that $\left( I^{10},T\right) $ is a spectral
pair in $\mathbb{R}^{10}$ such that $T-T$ does not contain one of the
canonical basis vectors; it follows that $T$ is not of the form (\ref{eq17}%
). In particular, we have verified
that the sets $\Lambda $
of the form (\ref{eq17}) do not suffice for cataloguing all
possible spectra for the unit cube $I^{10}$ in $\mathbb{R}^{10}$.

\section*{\label{Ap}Appendix: Spectral
Pairs of Measures}

\setcounter{section}{1}
\setcounter{theorem}{0}
\setcounter{equation}{0}
\renewcommand{\thesection}{\Alph{section}}%
We extend the concept of a spectral pair to a spectral pair of measures $%
(\mu ,\nu )$, where $\mu $ is a Borel measure on a locally compact abelian
group $G$ and $\nu $ a Borel measure on the dual group $\Gamma $. In the
past, we have mainly studied spectral pairs in the situation $G=\Gamma =%
\mathbb{R}^{d}$, $\left\langle x,\xi \right\rangle =e^{i2\pi x\xi }$, $\mu
\left( X\right) =m\left( \Omega \cap X\right) $, where $m$ denotes Lebesgue
measure, and $\Omega \subset \mathbb{R}^{d}$ is Lebesgue measurable with $%
m\left( \Omega \right) \ne 0$. In studying this situation, we found it
useful to also study spectral pairs in the cases $G=\mathbb{Z}^{d}$, $G=%
\mathbb{T}^{d}$, $G=\mathbb{Z}_{n}=\mathbb{Z}\diagup n\mathbb{Z}$, and $\mu $
a restriction of Haar measure on the respective groups $G$.

The present setup allows us to study all these cases simultaneously and to
expose a fundamental symmetry between the ``spectral set'' $\mu $ and the
``spectrum'' $\nu $ in a spectral pair $\left( \mu ,\nu \right) $ of
measures. We make the symmetric roles of $\mu $ and $\nu $ explicit, and we
show that $\mu \left( G\right) <\infty $ holds if and only if $\nu $ has an
atom.

Let $G$ be a locally compact abelian group (written additively). Let $m_{G}%
\colon G\rightarrow G$ be given by 
\begin{equation}
m_{G}\left( x\right) :=-x  \label{eqA.1}
\end{equation}
for $x\in G$. If $f\colon G\rightarrow \mathbb{C}$ is a complex-valued
function defined on $G$, then we let 
\begin{equation}
M_{G}f:=f\circ m_{G}.  \label{eqA.2}
\end{equation}

Let $\mu $ be a positive Borel measure on $G$.

Let $\tilde{\mu}:=\mu \circ m_{G}^{-1}$, i.e., 
\begin{equation}
\int f\,d\tilde{\mu}=\int M_{G}f\,d\mu .
\end{equation}
Then $M_{G}$ restricts to an isometric isomorphism of $\mathcal{L}^{2}\left(
\mu \right) $ onto $\mathcal{L}^{2}\left( \tilde{\mu}\right) $.

Let $\Gamma :=\hat{G}$ be the dual group of the group $G$, i.e., $\Gamma $
is the set of all continuous homomorphisms of $G$ into the unit circle $%
\mathbb{T}\simeq \mathbb{R}\diagup \mathbb{Z}$. Since $\hat{\Gamma}\simeq G$ 
\cite{HR63}, we can interpret $G$ as a set of homomorphisms on $\Gamma $. We
will write 
\begin{equation}
\left\langle x,\xi \right\rangle \in \mathbb{T},
\end{equation}
$x\in G$, $\xi \in \Gamma $ for the duality between $G$ and $\Gamma $. Then,
for each $\xi \in \Gamma $, $e_{\xi }\left( x\right) :=\left\langle x,\xi
\right\rangle $ determines a continuous homomorphism $G\rightarrow \mathbb{T}
$, and similarly $e_{x}\left( \xi \right) :=\left\langle x,\xi \right\rangle 
$ determines a continuous homomorphism $\Gamma \rightarrow \mathbb{T}$.

Define 
\begin{equation}
\left( Ff\right) \left( \xi \right) :=\int_{G}f\left( x\right) \overline{%
e_{\xi }\left( x\right) }\,d\mu \left( x\right)
\end{equation}
for $f\in \mathcal{L}^{1}\cap \mathcal{L}^{2}\left( \mu \right) $. Let $\nu $
be a second positive Borel measure on $\Gamma $. If $\{ Ff:f\in \mathcal{L}%
^{1}\cap \mathcal{L}^{2}\left( \mu \right) \} $ is a dense subset of $%
\mathcal{L}^{2}\left( \nu \right) $ and 
\begin{equation}
\int_{\Gamma}\left| Ff\right| ^{2}\,d\nu =\int_{G}\left| f\right| ^{2}\,d\mu
\end{equation}
for each $f\in \mathcal{L}^{1}\cap \mathcal{L}^{2}\left( \mu \right) $, then
we say $\left( \mu ,\nu \right) $ is a \emph{spectral pair} (\emph{of
measures\/}). In the affirmative case, $F$ ($=F_{\left( \mu ,\nu \right) }$)
extends, by continuity, to an isometric isomorphism of $\mathcal{L}%
^{2}\left( \mu \right) $ onto $\mathcal{L}^{2}\left( \nu \right) $.

Similarly to $m_{G}$ and $M_{G}$ in (\ref{eqA.1}) and (\ref{eqA.2}) above,
we introduce $m_{\Gamma }$ and $M_{\Gamma }$.

\begin{theorem}
\label{ThmA.1}If $\left( \mu ,\nu \right) $ is a spectral pair, then so is $%
\left( \tilde{\nu},\mu \right) $ with the transform $\tilde{\nu}$ from $\nu $
as introduced above.
\end{theorem}

\begin{proof}%
Let 
\begin{equation}
\left( \tilde{F}g\right) \left( x\right) =\int_{\Gamma }g\left( \xi \right) 
\overline{e_{x}\left( \xi \right) }\,d\tilde{\nu}\left( \xi \right) .
\end{equation}
We must show that $\tilde{F}$ extends, by continuity, to an isometric
isomorphism, mapping $\mathcal{L}^{2}\left( \tilde{\nu}\right) $ onto $%
\mathcal{L}^{2}\left( \mu \right) $. If $f\in \mathcal{L}^{1}\cap \mathcal{L}%
^{2}\left( \mu \right) $ and $g\in \mathcal{L}^{1}\cap \mathcal{L}^{2}\left( 
\tilde{\nu}\right) $, then 
\begin{equation}
\ip{M_{\Gamma }Ff}{g}%
_{\tilde{\nu}}
= 
\ip{f}{\tilde{F}g}%
_{\mu }
\end{equation}
by a simple
computation using
the fact
that $M_{\Gamma }$ is an
isometric isomorphism.
Hence, $\tilde{F}g=\left( M_{\Gamma }F\right) ^{*}g$. Since $M_{\Gamma }$
and $F$ both are isometric isomorphisms so is the adjoint $\left( M_{\Gamma
}F\right) ^{*}$.%
\end{proof}%

We have the following analogue of the usual Fourier inversion formula.

\begin{corollary}
\label{CorA.2}If $\left( \mu ,\nu \right) $ is a spectral pair of measures and if
$g\in \mathcal{L}^{1}\cap \mathcal{L}^{2}\left( \nu \right) 
$, then 
\begin{equation}
\left( F^{-1}g\right) \left( x\right) =\int_{\Gamma }g\left( \xi \right)
e_{x}\left( \xi \right) \,d\nu \left( \xi \right)   \label{eqAp.9}
\end{equation}
for $\mu $\textup{-a.e.}\ $x\in G$.
\end{corollary}

\begin{proof}%
In the proof of Theorem \ref{ThmA.1}, we showed that $\tilde{F}=\left(
M_{\Gamma }F\right) ^{*}$. Using $F^{-1}=F^{*}$ and $M_{\Gamma
}^{*}=M_{\Gamma }^{-1}=M_{\Gamma }^{{}}$ it follows that $F^{-1}=\tilde{F}%
M_{\Gamma }$,
and equation (\ref{eqAp.9}) is an immediate consequence.
\end{proof}%

The following result shows that every ``spectral set'' is a ``spectrum'' and
conversely that every ``spectrum'' is a ``spectral set''.

\begin{corollary}
\label{CorA.3}Let $\mu $ be a positive Borel measure on a locally compact
group $G$ and let $\nu $ be a positive Borel measure on the dual group $%
\Gamma =\hat{G}$. Then the following are equivalent:
\begin{enumerate}
\item  \label{CorA.3(1)}$\left( \mu ,\nu \right) $ is a spectral pair
on $\left( G,\Gamma \right) $,

\item  \label{CorA.3(2)}$\left( \tilde{\nu},\mu \right) $ is a spectral pair
on $\left( \Gamma ,G\right) $,

\item  \label{CorA.3(3)}$\left( \tilde{\mu},\tilde{\nu}\right) $ is a
spectral pair
on $\left( G,\Gamma \right) $,

\item  \label{CorA.3(4)}$\left( \nu ,\tilde{\mu}\right) $ is a spectral pair
on $\left( \Gamma ,G\right) $,
\end{enumerate}
where $\tilde{\mu}\left( \Delta \right) =\mu \left( -\Delta \right) $ and $%
\tilde{\nu}\left( \Delta ^{\prime }\right) =\nu \left( -\Delta ^{\prime
}\right) $, for any Borel sets $\Delta $ in $G$ and $\Delta ^{\prime }$ in $%
\Gamma $.
\end{corollary}


Corollary \ref{CorA.3} generalizes the well-known case where $G=\Gamma =%
\mathbb{R}^{d}$, $\mu =\nu =m$ and $F$ is the usual Fourier transform. For
example (\ref{CorA.3(3)}) corresponds to the fact that $\left( F^{2}f\right)
\left( \xi \right) =f\left( -\xi \right) $ for the usual Fourier transform.

\begin{theorem}
\label{th:symmetry}If $(\mu ,\nu )$ is a spectral pair, then so is $(\nu
,\mu )$.
\end{theorem}

\begin{proof}%
Suppose $(\mu ,\nu )$ is a spectral pair, then $(\nu ,\tilde{\mu })$ is
a spectral pair by \ref{CorA.3}, hence 
\begin{equation*}
F_{1}g(x):=\left\langle e_{x},g\right\rangle _{\nu }=\int g(\lambda )%
\overline{e_{x}(\lambda )}\,d\nu (\lambda )
\end{equation*}
determines an isometric isomorphism mapping $L^{2}(\nu )$ onto $L^{2}(%
\tilde{\mu })$. We must show that 
\begin{equation*}
F_{2}g(x):=\left\langle \overline{e_{x}},g\right\rangle _{\nu }=\int
g(\lambda )e_{x}(\lambda )\,d\nu (\lambda )
\end{equation*}
determines an isometric isomorphism mapping $L^{2}(\nu )$ onto $L^{2}(\mu )$
. But this is easy since 
\begin{eqnarray*}
F_{2}g(x) &=&\left\langle \overline{e_{x}},g\right\rangle _{\nu } \\
&=&\left\langle e_{-x},g\right\rangle _{\nu } \\
&=&F_{1}g(-x) \\
&=&M_{G}F_{1}g(x)
\end{eqnarray*}
so $F_{2}=M_{G}F_{1}$. By construction of $\tilde{\mu }$ from $\mu $ it
follows that $M_{G}$ is an isometric isomorphism mapping $L^{2}(\tilde{%
\mu })$ onto $L^{2}(\mu )$, hence $F_{2}=M_{G}F_{1}$, being the composition
of two isometric isomorphisms, is an isometric isomorphism as we needed to
show.%
\end{proof}%

If $\left( \mu ,\nu \right) $ is a spectral pair, then we may define a
strongly continuous unitary representation $U$ of $G$ on $\mathcal{L}%
^{2}\left( \mu \right) $ by 
\begin{equation}
\left( F\left( U\left( t\right) f\right) \right) \left( \xi \right)
:=e_{t}\left( \xi \right) \left( Ff\right) \left( \xi \right)
\end{equation}
for $f\in \mathcal{L}^{2}\left( \mu \right) $, $t\in G$ and $\nu $\textup{%
-a.e.}\ $\xi \in \Gamma $.

We have the following generalization of \cite[Corollary 1.11]{Ped87}, it
shows that if $\left( \mu ,\nu \right) $ is a spectral pair of measures and $%
\mu $ is a restriction of Haar measure to a set of finite measure, then the
pair $\left( \mu ,\nu \right) $ corresponds to a spectral pair of sets.

\begin{theorem}
\label{ThmA.4}Let $\left( \mu ,\nu \right) $ be a spectral pair. The
following are equivalent:

\begin{enumerate}
\item  \label{ThmA.4(1)}$\mu \left( G\right) <\infty $;

\item  \label{ThmA.4(2)}$\nu $ is a constant multiple of a counting measure;

\item  \label{ThmA.4(3)}$\nu \left( \left\{ \xi \right\} \right) \ne 0$ for
some $\xi \in \Gamma $.
\end{enumerate}

\noindent
In the affirmative case, the constant in \textup{(\ref{ThmA.4(2)})} is $\mu
\left( G\right) ^{-1}$.
\end{theorem}

\begin{proof}%
Since (\ref{ThmA.4(2)}) $\Rightarrow $ (\ref{ThmA.4(3)}) is trivial, we will
show that (\ref{ThmA.4(1)}) $\Rightarrow $ (\ref{ThmA.4(2)}) and (\ref
{ThmA.4(3)}) $\Rightarrow $ (\ref{ThmA.4(1)}).

\begin{proof}[Proof of \textup{(\ref{ThmA.4(1)}) $\Rightarrow $ (\ref{ThmA.4(2)})}]%
If $\mu \left( G\right) <\infty $, then $e_{\xi }\in \mathcal{L}^{2}\left(
\mu \right) $ for any $\xi \in \Gamma $. It follows that 
\begin{equation}
\ip{U\left( t\right) e_{\xi }}{e_{\eta }}%
_{\mu }= 
\ip{FU\left( t\right) e_{\xi }}{Fe_{\eta }}%
_{\mu }=\overline{e_{t}\left( \xi \right) }%
\ip{Fe_{\xi }}{Fe_{\eta }}%
_{\mu }=\overline{e_{t}\left( \xi \right) }%
\ip{e_{\xi }}{e_{\eta }}%
_{\mu }
\end{equation}
for $\nu $\textup{-a.e.}\ $\xi ,\eta \in \Gamma $. Now $U\left( t\right)
^{*}=U\left( -t\right) $ so $%
\ip{U\left( t\right) e_{\xi }}{e_{\eta }}%
_{\mu }=%
\ip{e_{\xi }}{U\left( -t\right) e_{\eta }}%
_{\mu }$ and therefore $\overline{e_{t}\left( \xi \right) }%
\ip{e_{\xi }}{e_{\eta }}%
_{\mu }=e_{-t}\left( \eta \right) 
\ip{e_{\xi }}{e_{\eta }}%
_{\mu }$ for any $t\in G$ and $\nu $\textup{-a.e.}\ $\xi ,\eta \in \Gamma $.
So, either $\xi =\eta $, or $%
\ip{e_{\xi }}{e_{\eta }}%
_{\mu }=0$. Consequently, 
\begin{equation}
\left( Fe_{\xi }\right) \left( \eta \right) =\int e_{\xi }\left( x\right) 
\overline{e_{\eta }\left( x\right) }\,d\mu \left( x\right) =%
\ip{e_{\eta }}{e_{\xi }}%
_{\mu }=\begin{cases} 0 & \text{ if } \eta \ne \xi \\ \mu \left( G\right) &
\text{ if } \eta =\xi . \end{cases}
\end{equation}
It follows that 
\begin{equation}
\mu \left( G\right) =%
\ip{e_{\xi }}{e_{\xi }}%
_{\mu }=%
\ip{Fe_{\xi }}{Fe_{\xi }}%
_{\nu }=\int \lvert Fe_{\xi }\left( \eta \right) \rvert^{2}\,d\nu \left(
\eta \right) =\mu \left( G\right) ^{2}\nu \left( \left\{ \xi \right\}
\right) .
\end{equation}
Hence, $\nu \left( \left\{ \xi \right\} \right) =\mu \left( G\right) ^{-1}$,
for any $\xi \in \limfunc{supp}\nu $.
\renewcommand{\qed}{}%
\end{proof}%

\begin{proof}[Proof of \textup{(\ref{ThmA.4(3)}) $\Rightarrow $ (\ref{ThmA.4(1)})}]%
By Theorem \ref{th:symmetry}
it is sufficient to show
that if $\mu\left( \left\{ x\right\} \right) \neq 0$
for some $x\in G$
then $\nu \left( 
\Gamma
\right) <\infty$.
Suppose first that
$x\in G$ is such that
$\mu \left( \left\{ x\right\} \right) \neq 1$. Since we can rescale
$\mu $ and $\nu $ by the same constant
we may assume $\mu \left( \left\{ x\right\} \right) =1$. Let
\begin{equation}
\delta _{x}\left( y\right) =\begin{cases} 1 & \text{ if } y=x \\ 0 & \text{
if } y\ne x. \end{cases}
\end{equation}
Then 
\begin{equation}
\left( F\delta _{x}\right) \left( \xi \right) =\int \delta _{x}\left(
y\right) \overline{e_{\xi }\left( y\right) }\,d\mu \left( y\right) =%
\overline{e_{\xi }\left( x\right) }=e_{m_{G}\left( x\right) }\left( \xi
\right) ,
\end{equation}
so $e_{m_{G}\left( x\right) }\in \mathcal{L}^{2}\left( \nu \right) $ and 
\begin{equation}
1=\left\| \delta _{x}\right\| _{\mu }^{2}=\left\| e_{m_{G}\left( x\right)
}\right\| _{\nu }^{2}=\nu \left( \Gamma \right) .
\end{equation}
An application of Theorem \ref{ThmA.1} completes the proof.%
\end{proof}%
\renewcommand{\qed}{}%
\end{proof}%

\begin{corollary}
\label{CorA.5}\raggedright
If $\left( \mu ,\nu \right) $ is a spectral pair such that $%
\mu \left( G\right) <\infty $, then $\left\{ e_{\xi }:\xi \in \limfunc{supp}%
\nu \right\} $ is an orthogonal basis for $\mathcal{L}^{2}\left( \mu \right) 
$.
\end{corollary}

\begin{proof}%
This is a simple consequence of the proof of Theorem \ref{ThmA.4}.%
\end{proof}%

Our next goal is to show that $\mu $ has uniform density. We first show that 
$U\left( t\right) $ acts by translation under appropriate circumstances.

\begin{lemma}
\label{LemA.6}Suppose $\left( \mu ,\nu \right) $ is a spectral pair. Let $%
\mathcal{O}\subset G$ be $\mu $-measurable and let $t\in G$. If $\mathcal{O}%
\cup \left( \mathcal{O}+t\right) \subset \limfunc{supp}\mu $, then 
\begin{equation}
\left( U\left( t\right) f\right) \left( x\right) =f\left( x+t\right)
\end{equation}
for $\mu $\textup{-a.e.}\ $x\in \mathcal{O}$ and every $f\in \mathcal{L}%
^{2}\left( \mu \right) $.
\end{lemma}

\begin{proof}%
If $Ff\in \mathcal{L}^{1}\cap \mathcal{L}^{2}\left( \nu \right) $, then 
\begin{align}
\left( U\left( t\right) f\right) \left( x\right) & =F^{-1}\left( e_{t}\left(
\xi \right) \left( Ff\right) \left( \xi \right) \right) \left( x\right) \\
& =\int e_{t}\left( \xi \right) e_{x}\left( \xi \right) \left( Ff\right)
\left( \xi \right) \,d\nu \left( \xi \right)  \notag \\
& =\left( F^{-1}Ff\right) \left( x+t\right)  \notag \\
& =f\left( x+t\right)  \notag
\end{align}
for $\mu $\textup{-a.e.}\ $x\in \mathcal{O}$.%
\end{proof}%

\begin{corollary}
\label{CorA.7}If $\left( \mu ,\nu \right) $ is a spectral pair, $t\in G$,
and $\mathcal{O}\subset G$ is $\mu $-measurable, then the inclusion $%
\mathcal{O}\cup \left( \mathcal{O}+t\right) \subset \limfunc{supp}\mu $
implies that 
\begin{equation}
\mu \left( \mathcal{O}\right) =\mu \left( \mathcal{O}+t\right) .
\end{equation}
\end{corollary}

\begin{proof}%
If $x\in \mathcal{O}+t$, then 
\begin{equation}
\left( U\left( -t\right) \chi _{\mathcal{O}}\right) \left( x\right) =\chi _{%
\mathcal{O}}\left( x-t\right) =1,
\end{equation}
hence 
\begin{equation}
\mu \left( \mathcal{O}\right) =\left\| \chi _{\mathcal{O}}\right\| _{\mu
}^{2}=\left\| U\left( -t\right) \chi _{\mathcal{O}}\right\| _{\mu }^{2}\le
\mu \left( \mathcal{O}+t\right) .
\end{equation}
Similarly, if $x\in \mathcal{O}$, then 
\begin{equation}
\left( U\left( t\right) \chi _{\mathcal{O}+t}\right) \left( x\right) =\chi _{%
\mathcal{O}+t}\left( x+t\right) =1,
\end{equation}
so 
\begin{equation}
\mu \left( \mathcal{O}+t\right) =\left\| \chi _{\mathcal{O}+t}\right\| _{\mu
}^{2}=\left\| U\left( t\right) \chi _{\mathcal{O}+t}\right\| _{\mu }^{2}\le
\mu \left( \mathcal{O}\right) .
\end{equation}
The desired equality is immediate.%
\end{proof}%

Our Corollary \ref{CorA.7} is related to the discussion in the recent paper 
\cite{KoLa96} of ``tiling the line by translates of a function'' as follows:
By Corollary \ref{CorA.7} such tilings do not come from spectral sets.

\begin{proposition}
\label{ProA.8}Let $\left( \mu ,\nu \right) $ be a spectral pair. Each $\xi
\in \Gamma $ is contained in an open set $\mathcal{O}_{\xi }$ with $\nu
\left( \mathcal{O}_{\xi }\right) <\infty $. In particular, each compact
subset of $\Gamma $ has finite $\nu $-measure; in short, $\nu $ is a Radon
measure.
\end{proposition}

\begin{proof}%
Let 
\begin{equation}
\mathcal{N}:=\left\{ \xi \in \Gamma :\left( Ff\right) \left( \xi \right) =0%
\text{ for all }f\in \mathcal{L}^{1}\cap \mathcal{L}^{2}\left( \mu \right)
\right\} .
\end{equation}
By density of $\mathcal{L}^{1}\cap \mathcal{L}^{2}\left( \mu \right) $ in $%
\mathcal{L}^{2}\left( \mu \right) $, we have $\nu \left( \mathcal{N}\right)
=0$. If $f\in \mathcal{L}^{1}\cap \mathcal{L}^{2}\left( \mu \right) $ then $%
Ff$ is continuous and 
\begin{align*}
\left( Ff\right) \left( \xi +\eta \right) & =\int f\left( x\right) \overline{%
e_{\xi +\eta }\left( x\right) }\,d\mu \left( x\right) \\
& =\int \left( \bar{e}_{\eta }f\right) \left( x\right) \overline{e_{\xi
}\left( x\right) }\,d\mu \left( x\right) \\
& =\left( F\left( \bar{e}_{\eta }f\right) \right) \left( \xi \right) ,
\end{align*}
so $\xi \in \mathcal{N}\Rightarrow \xi +\eta \in \mathcal{N}$, hence $%
\mathcal{N}=\varnothing $. Using $\mathcal{N}=\varnothing $, we see that
there exist $f_{\xi }\in \mathcal{L}^{1}\cap \mathcal{L}^{2}\left( \mu
\right) $ such that $\left( Ff_{\xi }\right) \left( \eta \right) \ge 1$ for $%
\eta $ in some neighborhood $\mathcal{O}_{\xi }$ of $\xi $ and $\nu \left( 
\mathcal{O}_{\xi }\right) \le \left\| Ff_{\xi }\right\| _{\nu }^{2}=\left\|
f_{\xi }\right\| _{\mu }^{2}$, as desired.%
\end{proof}%

\begin{proposition}
\label{ProA.9}Let $\left( \mu ,\nu \right) $ be a spectral pair. If $\mu
\left( G\right) <\infty $, then $\limfunc{supp}\nu $ is uniformly discrete
in the sense that there exists an open set $\mathcal{O}\subset \Gamma $ such
that $0\in \mathcal{O}$ and 
\begin{equation}
\left( \mathcal{O}+\xi \right) \cap \limfunc{supp}\nu =\left\{ \xi \right\}
\end{equation}
for all $\xi \in \limfunc{supp}\nu $.
\end{proposition}

\begin{proof}%
Let $g\left( \eta \right) :=%
\ip{e_{0}}{e_{\eta }}%
_{\mu }$. Then $g$ is continuous and $g\left( 0\right) =\mu \left( G\right)
^{\frac{1}{2}}$. Let $\mathcal{O}\subset \Gamma $ be an open set such that $%
0\in \mathcal{O}$ and $\lvert g\left( \eta \right) \rvert>\frac{1}{2}\mu
\left( G\right) ^{\frac{1}{2}}$ for any $\eta \in \mathcal{O}$. Then 
\begin{equation}
\lvert 
\ip{e_{\xi }}{e_{\eta }}%
\rvert=\lvert g\left( \eta -\xi \right) \rvert>\frac{1}{2}\mu \left(
G\right) ^{\frac{1}{2}}
\end{equation}
whenever $\eta \in \mathcal{O}+\xi $.%
\end{proof}%

\begin{theorem}[Uncertainty Principle]
\label{ThmA.U.P.}
Let $\left( \mu ,\nu \right) $ be a spectral pair
and $f\in L^{2}\left( \mu \right) $ with $\left\| f\right\| _{\mu }>0$. If
$A$, $B$ are measurable
sets in $G$ and $\Gamma $,
respectively, such that
\begin{align}
\left\| f-\chi _{A}f\right\| _{\mu }& \le \varepsilon \left\| f\right\|
_{\mu }\,,  \tag{i}  \label{ThmA.U.P.(1)} \\
\left\| Ff-\chi _{B}Ff\right\| _{\nu }& \le \delta \left\| Ff\right\| _{\nu
}=\delta \left\| f\right\| _{\mu }  \tag{ii}  \label{ThmA.U.P.(2)} \\
\TeXButton{both hold, then}{\intertext{both hold, then}}
\left( 1-\varepsilon -\delta \right) ^{2}
& \le
\mu \left( A\right) \nu \left( B\right) 
. 
\tag{iii}  \label{ThmA.U.P.(3)}
\end{align}
\end{theorem}

\begin{proof}[Proof \textup{(}sketch\textup{)}]%
If $f\in \mathcal{L}^{2}\left( \mu \right) $, $A\subset G$ and $B\subset
\Gamma $, then 
\begin{align}
\left\| \chi _{A}F^{-1}\chi _{B}Ff\right\| _{\mu }^{2}& =\int_{G}\chi
_{A}\left( x\right) \left| \int_{\Gamma }\chi _{B}\left( \xi \right)
e_{x}\left( \xi \right) \int_{G}f\left( y\right) \overline{e_{\xi }\left(
y\right) }\,d\mu \left( y\right) \,d\nu \left( \xi \right) \right|
^{2}\,d\mu \left( x\right) \\
& =\int \chi _{A}\left( x\right) \left| \int f\left( y\right) \int \chi
_{B}\left( \xi \right) e_{x-y}\left( \xi \right) \,d\nu \left( \xi \right)
\,d\mu \left( y\right) \right| ^{2}\,d\mu \left( x\right)  \notag \\
& =\int \chi _{A}\left( x\right) \left| \int f\left( y\right) \overline{%
\left( F^{-1}\chi _{B}\bar{e}_{x}\right) \left( y\right) }\,d\mu \left(
y\right) \right| ^{2}\,d\mu \left( x\right)  \notag \\
& \le \int \chi _{A}\left( x\right) \left\| f\right\| _{\mu }^{2}\left\|
F^{-1}\chi _{B}e_{x}\right\| _{\mu }^{2}\,d\mu \left( x\right)  \notag \\
& =\left\| f\right\| _{\mu }^{2}\int \chi _{A}\left( x\right) \left\| \chi
_{B}e_{x}\right\| _{\nu }^{2}\,d\mu \left( x\right)  \notag \\
& =\left\| f\right\| _{\mu }^{2}\mu \left( A\right) \nu \left( B\right) . 
\notag
\end{align}
If furthermore $\left\| f-\chi _{A}f\right\| _{\mu }\le \varepsilon \left\|
f\right\| _{\mu }$ and $\left\| Ff-\chi _{B}Ff\right\| _{\nu }\le \delta
\left\| Ff\right\| _{\nu }$, then also 
\begin{align}
\left\| f\right\| _{\mu }-\left\| \chi _{A}F^{-1}\chi _{B}Ff\right\| _{\mu
}& \le \left\| f-\chi _{A}F^{-1}\chi _{B}Ff\right\| _{\mu } \\
& \le \left\| f-\chi _{A}f\right\| _{\mu }+\left\| \chi _{A}f-\chi
_{A}F^{-1}\chi _{B}Ff\right\| _{\mu }  \notag \\
& \le \varepsilon \left\| f\right\| _{\mu }+\left\| f-F^{-1}\chi
_{B}Ff\right\| _{\mu }  \notag \\
& =\varepsilon \left\| f\right\| _{\mu }+\left\| Ff-\chi _{B}Ff\right\|
_{\nu }  \notag \\
& \le \varepsilon \left\| f\right\| _{\mu }+\delta \left\| f\right\| _{\mu }%
\text{,}  \notag \\
\intertext{so}%
\left( 1-\varepsilon -\delta \right) \left\| f\right\| _{\mu }& \le \left\|
\chi _{A}F^{-1}\chi _{B}Ff\right\| _{\mu }\,,  \notag
\end{align}
completing the proof.%
\end{proof}%

There is a vast literature on Uncertainty Principles. Theorem \ref{ThmA.U.P.}
and its proof are modelled after \cite{DoSt89}, see also \cite{deJe94}, \cite
{Be85}.
A comprehensive recent survey is \cite{FoSi97}.
A much more detailed analysis of operators similar to $\chi
_{A}F^{-1}\chi _{B}F$ appears in \cite{Lan67}.

\begin{corollary}
\label{CorA.U.P.}Let $\left( \mu ,\nu \right) $ be a spectral pair. If there
exists a sequence of $\mu $-measurable sets $A_{n}\subset G$ such that $\mu
\left( A_{n}\right) \ne 0$, and $\mu \left( A_{n}\right) \rightarrow 0$,
then $\nu \left( \Gamma \right) =+\infty $.
\end{corollary}

\begin{proof}%
Let $A\subset G$ be measurable with $0<\mu \left( A\right) <\infty $ and let 
$f=\mu \left( A\right) ^{-\frac{1}{2}}\chi _{A}$. Then $\left\| f\right\|
_{\mu }=1$ and 
\begin{equation}
1\le \mu \left( A\right) \nu \left( \Gamma \right)
\end{equation}
by the Uncertainty Principle with $B=\Gamma $. The desired conclusion is
immediate.%
\end{proof}%

Corollary \ref{CorA.U.P.} should be compared to
Theorem \ref{ThmA.4}. If $\left( \mu ,\nu \right) $ is a spectral
pair and $\mu $ is Lebesgue measure
restricted to a set $\Omega \subset \mathbb{R}^{d}$ of finite nonzero
Lebesgue measure, then
Theorem \ref{ThmA.4} and Corollary \ref{CorA.U.P.}
imply that $\nu $ is $\mu \left( \Omega \right) ^{-1}$ times
counting measure on an infinite
set $\Lambda \subset \mathbb{R}^{d}$. Finally, if $G$ is a
finite abelian group, $\Gamma $ is the
dual group,
$\mu $ is counting measure on $G$, and
$\nu $ is $\mu \left( G\right) ^{-1}$ times counting
measure on $\Gamma $, then  $\left( \mu ,\nu \right) $ is
a spectral pair; see \cite{HR63}. In particular,
the assumption in Corollary \ref{CorA.U.P.} that
there exist sets $A\subset G$ of arbitrarily
small $\mu $-measure cannot be removed.

\begin{acknowledgements}
The authors gratefully acknowledge excellent typesetting and graphics
production by Brian Treadway.
Very detailed
reports
from colleagues
guided us
towards substantial improvements in the exposition.
\end{acknowledgements}

\bibliographystyle{bftalpha}
\bibliography{jorgen}

\end{document}

%% file: inner.tex
\newbox\ipbox
\newcommand{\ip}[2]{\left\langle #1\mathrel{\mathchoice
{\setbox\ipbox=\hbox{$\displaystyle \left\langle\mathstrut #1#2\right\rangle$}
\vrule height\ht\ipbox width0.25pt depth\dp\ipbox}
{\setbox\ipbox=\hbox{$\textstyle \left\langle\mathstrut #1#2\right\rangle$}
\vrule height\ht\ipbox width0.25pt depth\dp\ipbox}
{\setbox\ipbox=\hbox{$\scriptstyle \left\langle\mathstrut #1#2\right\rangle$}
\vrule height\ht\ipbox width0.25pt depth\dp\ipbox}
{\setbox\ipbox=\hbox{$\scriptscriptstyle \left\langle\mathstrut #1#2\right\rangle$}
\vrule height\ht\ipbox width0.25pt depth\dp\ipbox}
} #2\right\rangle}
\newcommand{\diracb}[1]{\left\langle #1\mathrel{\mathchoice
{\setbox\ipbox=\hbox{$\displaystyle \left\langle\mathstrut #1\right.$}
\vrule height\ht\ipbox width0.25pt depth\dp\ipbox}
{\setbox\ipbox=\hbox{$\textstyle \left\langle\mathstrut #1\right.$}
\vrule height\ht\ipbox width0.25pt depth\dp\ipbox}
{\setbox\ipbox=\hbox{$\scriptstyle \left\langle\mathstrut #1\right.$}
\vrule height\ht\ipbox width0.25pt depth\dp\ipbox}
{\setbox\ipbox=\hbox{$\scriptscriptstyle \left\langle\mathstrut #1\right.$}
\vrule height\ht\ipbox width0.25pt depth\dp\ipbox}
}\right. }
\newcommand{\dirack}[1]{\left. \mathrel{\mathchoice
{\setbox\ipbox=\hbox{$\displaystyle \left.\mathstrut #1\right\rangle$}
\vrule height\ht\ipbox width0.25pt depth\dp\ipbox}
{\setbox\ipbox=\hbox{$\textstyle \left.\mathstrut #1\right\rangle$}
\vrule height\ht\ipbox width0.25pt depth\dp\ipbox}
{\setbox\ipbox=\hbox{$\scriptstyle \left.\mathstrut #1\right\rangle$}
\vrule height\ht\ipbox width0.25pt depth\dp\ipbox}
{\setbox\ipbox=\hbox{$\scriptscriptstyle \left.\mathstrut #1\right\rangle$}
\vrule height\ht\ipbox width0.25pt depth\dp\ipbox}
} #1\right\rangle}
\newcommand{\rip}[2]{\left( #1\mathrel{\mathchoice
{\setbox\ipbox=\hbox{$\displaystyle \left(\mathstrut #1#2\right)$}
\vrule height\ht\ipbox width0.25pt depth\dp\ipbox}
{\setbox\ipbox=\hbox{$\textstyle \left(\mathstrut #1#2\right)$}
\vrule height\ht\ipbox width0.25pt depth\dp\ipbox}
{\setbox\ipbox=\hbox{$\scriptstyle \left(\mathstrut #1#2\right)$}
\vrule height\ht\ipbox width0.25pt depth\dp\ipbox}
{\setbox\ipbox=\hbox{$\scriptscriptstyle \left(\mathstrut #1#2\right)$}
\vrule height\ht\ipbox width0.25pt depth\dp\ipbox}
} #2\right)}

%% file: tcibites.tex
\def\limfunc#1{\mathop{\rm #1}}%
\long\def\TeXButton#1#2{#2}%
\def\LaTeXparent#1{}%
\def\ChildStyles#1{}%